\documentclass{elsarticle}
\usepackage[utf8]{inputenc}
\usepackage[margin=1.5in]{geometry}
\usepackage{amsmath,amsthm,amssymb,amsfonts}
\usepackage{hyperref}

\DeclareMathOperator{\rank}{rank}

\newcommand{\R}{\mathbb{R}}
\newcommand{\E}{\mathbb{E}}
\newcommand{\norm}[1]{\left\lVert#1\right\rVert}

\newcounter{case}[section]
\newenvironment{case}[1][]
{\refstepcounter{case}\par\medskip\textbf{Case~\thecase. #1} \rmfamily}{\medskip}

\theoremstyle{plain}
 \newtheorem{thm}{Theorem}
 \newtheorem{conj}{Conjecture}
 \newtheorem{lem}{Lemma}
 \newtheorem*{cor}{Corollary}
\theoremstyle{remark}
 \newtheorem*{rem}{Remark}

\title{Few distance sets in $\ell_p$ spaces and $\ell_p$ product spaces}

\author[1]{Richard Chen}
\ead{rachen@mit.edu}
\author[2]{Feng Gui}
\ead{fenggui@mit.edu}
\author[3]{Jason Tang}
\ead{jtang21@belmontschools.net}
\author[4]{Nathan Xiong}
\ead{nxiong22@andover.edu}

\address[1]{Lexington High School, Lexington, MA 02421, USA}
\address[2]{Department of Mathematics, Massachusetts Institute of Technology, Cambridge, MA 02139, USA}
\address[3]{Belmont High School, Belmont, MA 02478, USA}
\address[4]{Phillips Academy Andover, Andover, MA 01810, USA}

\begin{document}

\begin{abstract}
    Kusner asked if $n+1$ points is the maximum number of points in $\R^n$ such that the $\ell_p$ distance $(1<p<\infty)$ between any two points is $1$. We present an improvement to the best known upper bound when $p$ is large in terms of $n$, as well as a generalization of the bound to $s$-distance sets. We also study equilateral sets in the $\ell_p$ sums of Euclidean spaces, deriving upper bounds on the size of an equilateral set for when $p=\infty$, $p$ is even, and for any $1\le p<\infty$.
\end{abstract}
\maketitle

\section{Introduction}
\subsection{Background}

A classic exercise in linear algebra asks for the maximum number of points in $\R^n$ such that the pairwise distances take only two values. One can associate a polynomial to each point in the set such that the polynomials are linearly independent. Then, one can show that the polynomials all lie in a subspace of dimension $(n+1)(n+4)/2$. Since the number of linearly independent vectors cannot exceed the dimension of the subspace, the cardinality of the set is at most $(n+1)(n+4)/2$. This beautiful argument was found by Larman, Rogers, and Seidel \cite{larman} in 1977, illustrating the power of algebraic techniques in combinatorics. 

We can ask the much more general question: ``In a metric space $X$, what is the maximum number of points such that the pairwise distances take only $s$ values?'' We use $e_s(X)$, or just $e(X)$ if $s=1$, to denote the answer to this question (by convention, we do not count $0$ as a distance). A set of points $S\subseteq X$ satisfying this question's conditions, i.e., the cardinality of the set $\{d(x,y):x,y\in S,x\neq y\}$ is $s$, is called an \emph{$s$-distance set}. A $1$-distance set is also known as an \emph{equilateral set}. Also, we typically restrict the metric space to a normed space, so that the specific distances used do not matter. Thus, we will always assume that the largest of the $s$ distances is $1$.

The posed question has been studied on many different spaces. The most famous result would be the upper bound $\binom{n+s}{s}$ on an $s$-distance set in $n$-dimensional Euclidean space, found by Bannai, Bannai, and Stanton \cite{bannai}. This result was also discovered independently by Blokhuis \cite{Blokhuis}. Another important case is when all points lie on the $n$-dimensional sphere. There is strong motivation for this problem as it has many applications in coding theory and design theory \cite{delsarte}. In particular, having tight upper bounds can help us find extremal configurations which satisfy unique properties. 

Similar results have been obtained in the hyperbolic space \cite{Blokhuis}, the Hamming space \cite{MusinNozaki}, and the Johnson space \cite{MusinNozaki}. Not as much is known in an arbitrary finite-dimensional normed space, known as a \emph{Minkowski space}, other than Petty's \cite{Petty}  general bound of $e(X)\le 2^n$, where $n=\dim X$. Swanepoel \cite{Swanepoel3} then conjectured that $e_s(X)\le (s+1)^n$ for a Minkowski space $X$ with dimension $n$ and proved it for the $n=2$ case. We should mention that equilateral sets in Minkowski spaces have been applied in differential geometry, where they are used to find minimal surfaces \cite{morgan}.

\subsection{Definitions}
In our paper, we investigate this problem on $\R^n$ with the $\ell_p$ norm, as well as on the $\ell_p$ sum of Euclidean spaces. For a point $x=(x_1,\ldots,x_n)\in \R^n$ and a $p\ge 1$, the $\ell_p$ norm is defined to be
\[\norm{x}_{p} = \left(\sum_{i=1}^n |x_i|^p\right)^\frac{1}{p},\]
and the $\ell_\infty$ norm is defined to be
\[\norm{x}_\infty = \max_{1\leq i\leq n} |x_i|.\]
The $\ell_1$ norm is the well-known ``taxicab'' norm, and the $\ell_2$ norm is the standard Euclidean norm. Throughout the paper, we write $\norm{\cdot}$ instead of $\norm{\cdot}_2$ to emphasize that we are using the Euclidean norm. We also use $\E^n$ to emphasize that we are in $n$-dimensional Euclidean space.
 
For Euclidean spaces $\E^{a_1},\ldots, \E^{a_n}$, we define their $\ell_p$ sum as the product space $\E^{a_1}\times \cdots \times \E^{a_n}$ equipped with the norm \[\norm{(x_1,\ldots, x_n)}_p=\left(\sum_{i=1}^n \norm{x_i}^p\right)^{\frac{1}{p}},\] where $x_i\in \E^{a_i}$ for $i=1,\ldots,n$. When $p=\infty$, the norm is
\[\norm{(x_1,\ldots, x_n)}_{\infty}=\max_{1\le i\le n}\norm{x_i}.\]

Then, for any $p\in[1,\infty]$, we let the distance between two points $x$ and $y$ be $\norm{x-y}_p$ and use $\E^{a_1}\oplus_p\cdots \oplus_p \E^{a_n}$ to describe this space.

Although our notation for the norm in $\ell_p$ spaces and in $\ell_p$ sums is the same, the norm being used should be clear from context.

\subsection{Previous Work and Our Results}

We first study $s$-distance sets in $\mathbb{R}^n$ equipped with the $\ell_p$ norm. This space is denoted by $\ell_p^n=(\R^n,\norm{\cdot}_p)$. The two most famous questions pertaining to this problem are Kusner's \cite{Guy} conjectures on equilateral sets.

\begin{conj}[Kusner]\label{conj:kusner1}
$e(\ell_1^n)=2n$.
\end{conj}
\begin{conj}[Kusner]\label{conj:kusner2}
$e(\ell_p^n)=n+1$ for $1<p<\infty$.
\end{conj}

For Conjecture~\ref{conj:kusner1}, note that the set $\{\pm e_i:i=1,\ldots,n\}$, where $e_i$ is the $i$-th standard basis vector, is equilateral in $\ell_1^n$, so $e(\ell_1^n)\ge 2n$. Currently, the best known upper bound is $e(\ell_1^n)\le cn\log n$ due to Alon and Pudl\'ak \cite{AP03}. It is also known that Conjecture~\ref{conj:kusner1} holds for $n=3$ (Bandelt, Chepoi, and Laurent \cite{bandelt}) and $n=4$ (Koolen, Laurent, and Schrijver \cite{koolen}).

As for Conjecture~\ref{conj:kusner2}, note that the set $\{e_1,\ldots,e_n,\lambda\sum_{i=1}^{n}e_i\}$ is equilateral for a suitable choice of $\lambda$, so $e(\ell_p^n)\ge n+1$. For $1<p<2$ and $n$ large enough, Swanepoel \cite{Swanepoel2} actually showed that $e(\ell_p^n)>n+1$, disproving Conjecture~\ref{conj:kusner2} in this case. The first nontrivial upper bound of $e(\ell_p^n)\le c_pn^{(p+1)/(p-1)}$ was found by Smyth \cite{smyth1} using an approximation argument. This result was later improved by Alon and Pudl\'ak \cite{AP03} and is currently the best known upper bound on $e(\ell_p^n)$ for arbitrary $n$ and $p$.
\begin{thm}[Alon and Pudl\'ak]\label{thm:alon}
For every $p\ge 1$, we have $e(\ell_p^n)\le c_pn^{(2p+2)/(2p-1)}$, where one may take $c_p=cp$ for an absolute $c>0$.
\end{thm}
Our first result is an improvement of Theorem~\ref{thm:alon} when $p$ is large in terms of $n$. One can check that $c>2$, so Theorem~\ref{thm:result1} is indeed an improvement. 
\begin{thm}\label{thm:result1}
Let $c>0$ be the constant from Theorem~\ref{thm:alon}. If $n>1$ and $p\ge c(n\log n)^2$, then $e(\ell_p^n)\le 2(p+1)n$.
\end{thm}
We should mention that when $p$ satisfies other special properties, Theorem~\ref{thm:alon} can be strengthened. If $p$ is an even integer, Swanepoel \cite{Swanepoel2} used a linear independence argument to show that
\[e(\ell_p^n)\le \begin{cases}
(\frac{p}{2}-1)n+1 & \text{if $p\equiv 0\pmod{4}$,}\\
\frac{p}{2}n+1&\text{if $p\equiv 2\pmod{4}$.}
\end{cases}\]
If $p$ is an odd integer, Alon and Pudl\'ak's argument for $e(\ell_1^n)$ extends to $e(\ell_p^n)$, giving the bound $e(\ell_p^n)<c_p n\log n$ in this case. 

Our next result is a generalization of Theorem~\ref{thm:alon} to $s$-distance sets. As far as we know, the only literature on $s$-distance sets in $\ell_p^n$ is by Smyth \cite{smyth2}. Our theorem below is strictly stronger than Conjecture 5 in \cite{smyth2}.
\begin{thm}\label{thm:result2}
If $s$ is a positive integer and $p$ is a real number satisfying $2p>s$, then $e_s(\ell_p^n)\le c_{p,s}n^{(2ps+2s)/(2p-s)}$ for a constant $c_{p,s}$ depending on $p$ and $s$.
\end{thm}
Our next three results are on equilateral sets in the $\ell_p$ sum of Euclidean spaces. As far as we know, this problem has not been well studied in this space. Swanepoel \cite{Swanepoel1} showed that $e(X\oplus_{\infty} Y)\le e(X)b_f(Y)$ for normed spaces $X$ and $Y$, where $b_f$ is the finite Borsuk number. However, explicit bounds are still unknown when $X$ and $Y$ are Euclidean spaces and when we take an $\ell_p$ sum instead of an $\ell_{\infty}$ sum. 
Our first result in this area almost completely resolves the problem for $\E^a\oplus_{\infty}\E^b$. Note that we have the obvious lower bound $e(\E^a\oplus_{\infty}\E^b)\ge e(\E^a)e(\E^b)= (a+1)(b+1)$ by taking the Cartesian product of the two equilateral sets. This lower bound is known to meet the upper bound when $a=2,3$ \cite{Swanepoel1}.
\begin{thm}\label{thm:result3}
Let $\E^a$ and $\E^b$ be Euclidean spaces. Then, $e(\E^a\oplus_{\infty}\E^b)\le (a+1)(b+1)+1$.
\end{thm}
Our second result in this area provides an upper bound when $p$ is even. 
\begin{thm}\label{thm:result4}
Let $\E^a$ and $\E^b$ be Euclidean spaces, and let $p$ be an even integer. Then, $e(\E^a\oplus_p\E^b)\le \binom{a+p/2}{a}+\binom{b+p/2}{b}$.
\end{thm}
Finally, we extend Alon and Pudl\'ak's \cite{AP03} result on equilateral sets in $\ell_p^n$ to certain $\ell_p$ product spaces. We present an upper bound on the $\ell_p$ sum of $n$ Euclidean spaces for any $1\le p<\infty$. Observe that Theorem~\ref{thm:alon} is a special case of our theorem below when $a_1=\cdots = a_n=1$.
\begin{thm}\label{thm:result5}
Let $\E^{a_1},\ldots, \E^{a_n}$ be Euclidean spaces and set $\displaystyle a=\max_{1\le i\le n}a_i$. If $2p>a$, then $e(\mathbb{E}^{a_1}\oplus_p \cdots \oplus_p \mathbb{E}^{a_n})\le c_{p,a}n^{\frac{2p+2a}{2p-a}}$ for a constant $c_{p,a}$ depending on $p$ and $a$.
\end{thm}

Our paper is structured as follows. In Section~\ref{sec:prelim}, we introduce two important tools used in our proofs. Then, we prove Theorems~\ref{thm:result1}, \ref{thm:result2}, \ref{thm:result3}, \ref{thm:result4}, and  \ref{thm:result5} in Sections~\ref{sec:proof1}, \ref{sec:proof2}, \ref{sec:proof3}, \ref{sec:proof4}, and \ref{sec:proof5} respectively.

\section{Preliminaries}\label{sec:prelim}

In this section, we present two famous results that we use later. The first is the Rank Lemma, which allows us to estimate the rank of a symmetric matrix. The second is Jackson's Theorem, a celebrated result from approximation theory. This theorem allows us to approximate $|x|^p$ by a polynomial with sufficiently small error.

\subsection{The Rank Lemma}
\begin{lem}[Rank Lemma]\label{lem:ranklemma}
For a real symmetric $n\times n$ nonzero matrix $A$,
\[\rank A \ge \frac{(\sum_{i=1}^{n} a_{ii})^2}{\sum_{i,j=1}^n a_{ij}^2}.\]
\end{lem}
We will frequently make use of the following corollary.
\begin{cor}\label{cor:ranklemma}
Let $A$ be a real symmetric $n\times n$ matrix with $a_{ii}=1$ and $|a_{ij}|\le \varepsilon$ for all $i\neq j$. Then \[\rank A\ge \frac{n}{1+(n-1)\varepsilon^2}.\]
Choosing $\varepsilon=n^{-1/2}$ gives $\rank A\ge n/2$.
\end{cor}

\subsection{Jackson's Theorem}
\begin{thm}[Jackson,  \cite{jackson}]\label{thm:Jackson}
For any $f\in C^k[-1,1]$ and positive integer $d$, there exists a polynomial $P$ with degree at most $d$ such that
\[|f(x)-P(x)|\le \frac{c^k}{(d+1)_k}\omega(f^{(k)},1/d) \qquad\text{for all }x\in [-1,1],\]
where $c>0$ is an absolute constant, $\omega(f,\delta):=\sup\{|f(x)-f(y)|:|x-y|\le \delta\}$ denotes the modulus of continuity of $f$, and $(d+1)_k=(d+1)d\cdots (d-k+2)$ uses falling factorial notation.
\end{thm}
From Theorem~\ref{thm:Jackson}, we can recover the following lemma. This lemma was first given and proved in \cite{Smyth_thesis} and later also transmitted in \cite{MR2117069}.
\begin{lem}\label{lem:jackson}
For any $p\ge 1$ and $d\ge \lceil p\rceil$, there exists a polynomial $P$ with degree at most $d$ such that
\begin{equation}\label{eq:jackson}
|P(x)-|x|^p|\le \frac{B(p)}{d^p}\qquad\text{for all }x\in [-1,1]
\end{equation}
where $B(p)=(\lceil p\rceil^p(1+\pi^2/2)^{\lceil p\rceil}(p)_{\lceil p\rceil-1})/\lceil p\rceil!$.
\end{lem}
We will always assume that the polynomial $P$ in Lemma~\ref{lem:jackson} is even and that $P(0)=0$. If $P$ is not even, we can take the even part of $P$. If $P(0)\neq 0$, we can take the polynomial $Q(x)=P(x)-P(0)$ as this only increases the error term by a factor of $2$.

\section{Bound on equilateral sets in $\ell_p^n$ for large $p$}\label{sec:proof1}
We start with an important lemma about the $\ell_p$ norm.
\begin{lem}\label{lem:holderbound}
Suppose $1\le p\le q$. Then for any $x\in \R^n$,
\[\norm{x}_q\le \norm{x}_p\le n^{\frac{1}{p}-\frac{1}{q}}\norm{x}_q.\]
\end{lem}
\begin{proof}
The left hand inequality is a well-known property of the $\ell_p$ norm. The right hand inequality is the Power Mean Inequality, which can be proven using Jensen's Inequality.
\end{proof}

Suppose our equilateral set is $\{p_1,\ldots,p_m\}$. Let $k$ be the closest even integer to $p$, rounding up if $p$ is odd. The main idea is that by Lemma~\ref{lem:holderbound}, we can approximate $\norm{\cdot}_p$ by $\norm{\cdot }_k$. If $p$ is large enough in terms of $n$, the error term is sufficiently small. From there, it suffices to bound the size of an approximately-equilateral set in $\ell_k^n$, which can be done with a linear algebra argument.

\begin{proof}[Proof of Theorem~\ref{thm:result1}]
We use the notation explained above. There are two cases to consider depending on whether $k$ is greater than or less than $p$.

\begin{case}\label{case:holder1}
We have $\lfloor p\rfloor $ is odd, so $p<k$.
\end{case}

From our bound on $p$ and the fact that $p\ge 2$,
\begin{align*}
    m &\le cpn^{(2p+2)/(2p-1)}\\
    &\le cpn^2\\
    &\le \frac{p^2}{(\log n)^2}\\
    &\le (1-n^{-1/p})^{-2},
\end{align*}
where the first inequality holds by Theorem~\ref{thm:alon} and the last inequality by the inequality $1-1/x\le \log(x)$ valid for $x>0$. Assuming without loss of generality that $m>1$, we may rearrange the result above into
\[-p\log_n\left(1-\frac{1}{\sqrt{m}}\right)\ge 1.\]
Because $k$ is the smallest even integer greater than $p$, we have
\begin{equation}\label{eq:kbound1}
p< k\le p+1\le p-p\log_n\left(1-\frac{1}{\sqrt{m}}\right).
\end{equation}
Now, embed the $m$ points into $\ell_k^n$. Since $p<k$, Lemma~\ref{lem:holderbound} implies
\begin{equation}\label{eq:holderbound1}
 n^{\frac{1}{k}-\frac{1}{p}}\le \norm{p_i-p_j}_k\le 1
\end{equation}
for all $i\neq j$. Consider the $m$ functions $f_i:\R^n\to \R$ given by
$f_i(x)=1-\norm{p_i-x}_k^k$, and
let $A$ be the matrix with $a_{ij}=f_i(p_j)$. Clearly, $a_{ii}=1$, and from (\ref{eq:kbound1}) and (\ref{eq:holderbound1}),
\begin{align*}
    |a_{ij}|&\le
    1-n^{1-\frac{k}{p}}\\
    &\le 1-n^{1-(1-\log_n(1-1/{\sqrt{m}}))}\\
    &= \frac{1}{\sqrt{m}},
\end{align*}
for all $i\neq j$. Since $A$ is symmetric, Lemma~\ref{lem:ranklemma} tells us that $\rank A\ge m/2$. For the upper bound on the rank, note that every $f_i$ lies in the span of the set of polynomials
\[\mathcal{S}=\{1,x_1,\ldots,x_n,x_1^2,\ldots,x_n^2,\ldots,x_1^{k-1},\ldots,x_n^{k-1},\sum_{t=1}^{n}x_t^k\}.\]
The $i$-th row vector of $A$ is $(f_i(p_1),\ldots,f_i(p_m))$. Hence, every row vector belongs to the subspace spanned by $\{(f(p_1),\ldots,f(p_m)):f\in \mathcal{S}\}$, which has dimension at most $|\mathcal{S}|=(k-1)n+2\le kn$. It follows that $\rank A\le kn\le (p+1)n$. Combining the upper and lower bound, we have
\[m\le 2(p+1)n,\]
as desired.

\begin{case}\label{case:holder2}
We have $\lfloor p\rfloor $ is even, so $p\ge k$.
\end{case}

Since $c,n\ge 2$, $p>4$. So, similar to Case~\ref{case:holder1}, we have
\begin{align*}
    m&\le cpn^{(2p+2)/(2p-1)}\\
    &\le cpn^{2-2/p}\\
    &\le n^{-2/p}\cdot\frac{p^2}{(\log n)^2}\\
    &\le n^{-2/p}(1-n^{-1/p})^{-2}\\
    &= (n^{1/p}-1)^{-2}.
\end{align*}
Rearranging the result gives us
\[p\log_n\left(1+\frac{1}{\sqrt{m}}\right)\ge 1.\]
Because $k$ is the largest even integer less than or equal to $p$, we have
\begin{equation}\label{eq:kbound2}
p\ge k\ge p-1\ge p-p\log_n\left(1+\frac{1}{\sqrt{m}}\right).
\end{equation}
Embed the $m$ points into $\ell_k^n$. This time, Lemma~\ref{lem:holderbound} implies
\begin{equation}\label{eq:holderbound2}
1\le \norm{p_i-p_j}_k\le n^{\frac{1}{k}-\frac{1}{p}}.
\end{equation}
We define $f_i$ and $A$ like in Case~\ref{case:holder1}. Again, $a_{ii}=1$, and from (\ref{eq:kbound2}) and (\ref{eq:holderbound2}),
\begin{align*}
    |a_{ij}|&\le n^{1-\frac{k}{p}}-1\\
    &\le n^{1-(1-\log_n(1+1/{\sqrt{m}}))}-1\\
    &= \frac{1}{\sqrt{m}},
\end{align*}
for all $i\neq j$. Applying Lemma \ref{lem:ranklemma}, $\rank A\ge m/2$. So, similar to Case~\ref{case:holder1}, we have the bound $\rank A\le kn\le pn$. The final result is
\[m\le 2pn<2(p+1)n.\]
Having considered all cases, the proof is complete.
\end{proof}

\section{Bound on $s$-distance sets in $\ell_p^n$}\label{sec:proof2}
Smyth \cite{Smyth_thesis}, see also \cite{smyth2}, first combined Jackson's theorem with the invertibility of diagonally dominated matrices to prove $e(\ell_p^n)\le cn^{(p+1)/(p-1)}$ for $p>1$. Alon and Pudl\'ak \cite{AP03} then improved this bound by substituting the rank lemma for the invertibility lemma to achieve Theorem~\ref{thm:alon}. We would like to extend this method of proof to the $s$-distance case. However, as pointed out by Smyth \cite{smyth2}, if one of the distances is arbitrarily small, we need arbitrarily high degrees of approximation. So, we want to impose a lower bound on our distances. 

Consider a two-distance set with distances $1$ and $a$, where $a$ is very small. Intuitively, this set should look like several ``clusters'' of points, such that the distance between each cluster is $1$. So, if we look globally, this set looks like an equilateral set with distance $1$, but if we look locally at each cluster, we have an equilateral set with distance $a$. This means that we have essentially reduced the problem from one about two-distance sets to one about equilateral sets. We carry this intuition to $s$-distance sets and formalize this argument with an induction.

\begin{proof}[Proof of Theorem~\ref{thm:result2}]
We use strong induction on $s$. The base case is $s=1$, which was proven by Alon and Pudl\'ak \cite{AP03}. For the inductive step, assume that the statement holds true for all $1\le s<k$. We prove that it is true for $s=k$.

Let $S$ be a $k$-distance set in $\ell_p^n$. Let our $k$ distances be $1=a_1>a_2>\cdots>a_k$. There are two cases to consider depending on whether the $k$ distances are lower bounded or not.

\begin{case}\label{case:distance1}
The smallest distance $a_k$ is less than $2^{1-k}$.
\end{case}

There exists an index $1\le i< k$ such that $a_{i}>2a_{i+1}$. Let  $S'\subseteq S$ be a maximal $i$-distance set using the distances $a_1,\ldots,a_i$. Every point $p\in S\setminus S'$ is a distance at most $a_{i+1}$ from some point in $S'$, otherwise $p$ can be added to $S'$, contrary to the maximality of $S'$. So, if we draw a closed ball of radius $a_{i+1}$ around every point in $S'$, every point in $S$ lies within some ball. These balls are disjoint from the condition $a_i>2a_{i+1}$. Furthermore, by the triangle inequality, within each ball, the distance between any two points is at most $2a_{i+1}<a_i$, and thus is at most $a_{i+1}$.  It follows that within every ball, we have a $k'$-distance set, where $k'\le k-i$, using $k'$ of the distances from $a_{i+1},\ldots,a_{k}$. This implies the bound
\[e_k(\ell_p^n) \le e_i(\ell_p^n)\cdot\max_{0\le j\le k-i}e_j(\ell_p^n).\]
Applying the inductive hypothesis gives us
\begin{align*}
    e_k(\ell_p^n) &\le c_{p,i}n^{\frac{2pi+2i}{2p-i}}\cdot \max_{0\le j\le k-i}c_{p,j}n^{\frac{2pj+2j}{2p-j}}\\
    &\le c_{p,i}n^{\frac{2pi+2i}{2p-i}}\cdot c_{p,k-i}n^{\frac{2p(k-i)+2(k-i)}{2p-(k-i)}}\\
    &\le c_{p,k}n^{\frac{2pk+2k}{2p-k}},
\end{align*}
as desired. The first inequality follows by analyzing the function
\[f_1(x)=\frac{2px+2x}{2p-x}\]
on the interval $[0,2p)$. Its first derivative is
\[f_1'(x)=\frac{4p(p+1)}{(2p-x)^2},\]so $f_1$ is increasing on $[0,2p)$ and hence on $[0,k-i]$. The second inequality follows by analyzing the related function 
\[f_2(x)=\frac{2px+2x}{2p-x}+\frac{2p(k-x)+2(k-x)}{2p-(k-x)}\]
on the interval $[0,k]$. It is symmetric about $x=\tfrac{k}{2}$, and its first derivative is
\[f_2'(x)=\frac{4p(p+1)}{(2p-x)^2}-\frac{4p(p+1)}{(2p-(k-x))^2}.\]
Since $2p>k$, $f_2$ is decreasing on $[0,\tfrac{k}{2}]$ and increasing on $[\tfrac{k}{2},k]$. Thus, it is maximised at $x=0$ and $x=k$. 

\begin{case}\label{case:distance2}
The smallest distance $a_k$ is at least $2^{1-k}$.
\end{case}

Suppose $S$ consists of the points $\{p_1,\ldots,p_m\}$. For convenience, let $\pi=a_1^pa_2^p\cdots a_k^p$. Fix $B(p)$ as the constant from Lemma~\ref{lem:jackson}. Define $c$ as
\[c = \max(B(p)k\cdot 2^{pk^2-pk+2k},(2^{1/p}-1)^{-p}).\]
Then, let $d$ be a positive integer satisfying
\[cn\sqrt{m}<d^p<2cn\sqrt{m},\]
which is possible since $c\ge (2^{1/p}-1)^{-p}$
.

Lemma~\ref{lem:jackson} allows us to pick an even polynomial $P$ with $P(0)=0$ and degree at most $d$ such that \[|P(x)-|x|^p|\le \frac{B(p)}{d^p}\qquad\text{for all }x\in [-1,1].\] Consider the $m$ functions $f_i:\R^n\to \R$ given by
\[f_i(x)=\frac{1}{\pi}\prod_{u=1}^{k}\left(a_u^p-\sum_{t=1}^{n}|x_t-p_{it}|^p\right),\]
and their polynomial approximations
\[g_i(x)=\frac{1}{\pi}\prod_{u=1}^{k}\left(a_u^p-\sum_{t=1}^{n}P(x_t-p_{it})\right).\]
Let $A$ be the $m\times m$ matrix given by $a_{ij}=g_i(p_j)$. First, since $P(0)=0$, $a_{ii}=1$ for all $i$. We now estimate $a_{ij}$ for $i\neq j$. Expand \[f_i(x)=\frac{1}{\pi}\sum_{\ell=0}^{k}(-1)^{k+\ell}\sigma(\ell)\left(\sum_{t=1}^{n}|x_t-p_{it}|^p\right)^{k-\ell},\]
and 
\[g_i(x)=\frac{1}{\pi}\sum_{\ell=0}^{k}(-1)^{k+\ell}\sigma(\ell)\left(\sum_{t=1}^{n}P(x_t-p_{it})\right)^{k-\ell},\]
where $\sigma(\ell)$ denotes the $\ell$-th elementary symmetric polynomial in $a_1^p,\ldots,a_k^p$. For convenience, we define $X_{ij}=\sum_{t=1}^{n}|p_{jt}-p_{it}|^p$ and $Y_{ij}=\sum_{t=1}^{n}P(p_{jt}-p_{it})$. Applying the triangle inequality with (\ref{eq:jackson}), we have
\[|X_{ij}-Y_{ij}|\le \frac{nB(p)}{d^p}.\]
Since $i\neq j$, recall that $X_{ij}\ge 2^{(1-k)p}$. Combining this with the fact that $c>B(p)\cdot 2^{(k-1)p}$ implies that $Y_{ij}$ is positive. Now, we are ready to upper bound $|a_{ij}|$. Since $f_i(p_j)=0$ for $i\neq j$, we have $a_{ij}=g_i(p_j)=g_i(p_j)-f_i(p_j)$. Thus
\begin{align*}
    |a_{ij}| &= \frac{1}{\pi}\left |\sum_{\ell=0}^{k}(-1)^{k+\ell}\sigma(\ell)(Y_{ij}^{k-\ell}-X_{ij}^{k-\ell})\right |\\
    &\le \frac{1}{\pi}\cdot\frac{nB(p)}{d^p}\sum_{\ell=0}^{k-1}\sigma(\ell)\sum_{r=1}^{k-\ell}|X_{ij}^{k-\ell-r}Y_{ij}^{r-1}|\\
    &\le \frac{1}{\pi}\cdot\frac{nB(p)}{d^p}\sum_{\ell=0}^{k-1}\sigma(\ell)\sum_{r=1}^{k-\ell}\left |\frac{(k-\ell-r)X_{ij}+(r-1)Y_{ij}}{k-\ell-1}\right |^{k-\ell-1}\\
    &= \frac{1}{\pi}\cdot\frac{nB(p)}{d^p}\sum_{\ell=0}^{k-1}\sigma(\ell)\sum_{r=1}^{k-\ell}\left |\frac{(r-1)(Y_{ij}-X_{ij})}{k-\ell-1}+X_{ij}\right |^{k-\ell-1}\\
    &\le \frac{1}{\pi}\cdot\frac{nB(p)}{d^p}\sum_{\ell=0}^{k-1}\sigma(\ell)(k-\ell)\left(\frac{nB(p)}{d^p}+1\right)^{k-\ell-1}\\
    &\le \frac{1}{\pi}\cdot\frac{nB(p)k}{d^p}\cdot 2^k\sum_{\ell=0}^{k-1}\sigma(\ell).
\end{align*}
The first inequality is achieved by the triangle inequality, the factorization $x^r-y^r=(x-y)\sum_{i=0}^{r-1}x^{r-i}y^i$, and the upper bound on $|X_{ij}-Y_{ij}|$. The second inequality follows by the AM-GM inequality. The last inequality holds because $nB(p)/d^p<B(p)/c<1$ by our choice of $c$ and $d$, whence $k2^k\ge (k-l)(1+nB(p)/d^p)^{k-\ell-1}$.

Now, since $a_u\ge 2^{1-k}$, we can lower bound $\pi$ with
\[\pi=a_1^pa_2^p\cdots a_k^p\ge (2^{(1-k)p})^{k}= 2^{pk-pk^2}.\]
On the other hand, since $a_t\le 1$, we can upper bound $\sigma(\ell)$ with
\[\sigma(\ell)=\sum_{1\le j_1<j_2<\cdots<j_{\ell}\le k}a_{j_1}^pa_{j_2}^p\cdots a_{j_{\ell}}^p\le \binom{k}{\ell}.\]
This gives us $\sum_{\ell=0}^{k-1}\sigma(\ell)< 2^{k}$. Putting everything together, we have
\begin{align*}
    |a_{ij}| &< 2^{pk^2-pk+2k}\cdot\frac{nB(p)k}{d^p}\\
    &< \frac{1}{\sqrt{m}}
\end{align*}
from our choice of $c$ and $d$.

Recall that $P$ is even, so the matrix $A$ is symmetric. Lemma~\ref{lem:ranklemma} then tells us that $\rank A\ge m/2$. We now find an upper bound for the rank.

Note that the polynomials $a_u^p-\sum_{t=1}^{n}P(x_t-p_{it})$ belong to the span of the set
\[\{1,x_1,\ldots,x_n,x_1^2,\ldots,x_n^2,\ldots,x_1^{d-1},\ldots,x_n^{d-1},\sum_{t=1}^{n}x_t^d\}.\]
This set consists of $(d-1)n+2\le dn$ polynomials. Thus, all the $g_i$ belong to the span of a set $\mathcal{S}$ of at most $(dn)^k$ polynomials. The $i$-th row vector of $A$ is $(g_i(p_1),\ldots,g_i(p_m))$. Hence, every row vector belongs to the subspace spanned by
\[\{(f(p_1),\ldots,f(p_m)):f\in \mathcal{S}\},\]
which has dimension at most $|\mathcal{S}|\le (dn)^k$. This implies $\rank A\le (dn)^k$.

The upper and lower bounds combine to give $(dn)^k\ge m/2$. Using the upper bound on $d^p$ and the condition $2p>k$, we can rearrange the inequality to obtain
\[m\le c_{p,k}n^{(2pk+2k)/(2p-k)}.\]
This completes the induction.
\end{proof}
 
\section{Bound on equilateral sets in $\ell_\infty$ product spaces}\label{sec:proof3}

\begin{proof}[Proof of Theorem~\ref{thm:result3}] Write $x=(\widetilde{x}_1,\widetilde{x}_2)$ for each point $x\in \E^a\oplus_{\infty}\E^b$, where $\widetilde{x}_1\in \E^a$ and $\widetilde{x}_2\in \E^b$. Let $S$ be our equilateral set with cardinality $m$. Consider the $m$ functions $f_u:\R^{a+b}\to \R$ defined by
\[ f_u(x)=\left(1-\norm{\widetilde{x}_1-\widetilde{u}_1}^2\right)\left(1-\norm{\widetilde{x}_2-\widetilde{u}_2}^2\right),\]
for all $u\in S$. Note that $f_u(v)=\delta_{uv}$ for all $u,v \in S$, so the $f_u$ are linearly independent.

We can expand $f_u$ as
\begin{align*}
    f_u(x) &= \left(1-\sum_{t=1}^{a}(\widetilde{x}_{1t}-\widetilde{u}_{1t})^2\right)\left(1-\sum_{t=1}^{b}(\widetilde{x}_{2t}-\widetilde{u}_{2t})^2\right)\\
    &= \left(1-\norm{\widetilde{u}_1}^2-\norm{\widetilde{x}_1}^2+2\sum_{t=1}^{a}\widetilde{x}_{1t}\widetilde{u}_{1t}\right)\left(1-\norm{\widetilde{u}_2}^2-\norm{\widetilde{x}_2}^2+2\sum_{t=1}^{b}\widetilde{x}_{2t}\widetilde{u}_{2t}\right)
\end{align*}
So, the $f_u$ are all spanned by the following set of $(a+2)(b+2)$ polynomials \[\{1,\widetilde{x}_{1i},\widetilde{x}_{2j},\widetilde{x}_{1i}\widetilde{x}_{2j},\norm{\widetilde{x}_1}^2,\norm{\widetilde{x}_2}^2, \norm{\widetilde{x}_2}^2\widetilde{x}_{1i},  \norm{\widetilde{x}_1}^2\widetilde{x}_{2j},\norm{\widetilde{x}_1}^2\norm{\widetilde{x}_2}^2: 1\le i\le a, 1\le j\le b\}.\]
This implies the bound $m\le (a+2)(b+2)$.

We will now prove that the set of polynomials $\{f_u,1,x_k,\norm{\widetilde{x}_1}^2:u\in S, 1\le k\le a+b\}$ is  linearly independent. Assume for the sake of contradiction that we have a dependence
\begin{equation}\label{eq:dependence}
\sum_{u\in S}\alpha_uf_u+\sum_{k=1}^{a+b}\beta_kx_k+\gamma\norm{\widetilde{x}_1}^2+\delta=0.
\end{equation}
The left hand side of (\ref{eq:dependence}) is a polynomial in the $x_k$ that is identically zero. Thus, extracting the coefficient of $\norm{\widetilde{x}_1}^2\norm{\widetilde{x}_2}^2$, we have
\begin{equation}\label{eq:extract1}
    \sum_{u\in S}\alpha_u=0.
\end{equation}
Extracting the coefficient of $\norm{\widetilde{x}_r}^2x_k$, for the appropriate $r\in\{1,2\}$, we have
\begin{equation}\label{eq:extract2}
    \sum_{u\in S}\alpha_uu_{k}=0.
\end{equation}
Extracting the coefficient of $\norm{\widetilde{x}_2}^2$ and applying (\ref{eq:extract1}), we have 
\begin{equation}\label{eq:extract3}
    \sum_{u\in S}\alpha_u\norm{\widetilde{u}_1}^2=0.
\end{equation}
Now, plug $u$ into (\ref{eq:dependence}), multiply both sides by $\alpha_u$, and sum over all $u\in S$.
\[ \sum_{u\in S}\alpha_u^2+\sum_{k=1}^{a+b}\beta_k\sum_{u\in S}\alpha_uu_{k}+\gamma\sum_{u\in S}\alpha_u\norm{\widetilde{u}_1}^2+\delta\sum_{u\in S}\alpha_u=0.\]
Applying (\ref{eq:extract1}), (\ref{eq:extract2}), and (\ref{eq:extract3}) implies $\alpha_u=0$ for all $u\in S$. It easily follows that all other coefficients are zero, as desired.

Now, we know that $m+a+b+2\le (a+2)(b+2)$. This rearranges into $m\le (a+1)(b+1)+1$.
\end{proof}

\section{Bound on equilateral sets in $\ell_p$ product spaces for even $p$}\label{sec:proof4}
\begin{proof}[Proof of Theorem~\ref{thm:result4}] Let $S$ be an equilateral set in $\E^a \oplus_p \E^b$ with $m$ points. For every $u\in S$, define the function $f_u:\R^{a+b}\to \R$ by
\[f_u(x)=1-\norm{\widetilde{x}_1-\widetilde{u}_1}^p-\norm{\widetilde{x}_2-\widetilde{u}_2}^p\]
for all $x=(\widetilde{x_1},\widetilde{x_2})\in \E^a\oplus_p \E^b$ so that $f_u(v)=\delta_{uv}$ for all $u,v\in S$. It follows that the $m$ polynomials are linearly independent. Now, we can expand $f_u$ as
\begin{align*}
    f_u(x) &= 1-\left(\sum_{t=1}^{a}(\widetilde{x}_{1t}-\widetilde{u}_{1t})^2\right)^{p/2}-\left(\sum_{t=1}^{b}(\widetilde{x}_{2t}-\widetilde{u}_{2t})^2\right)^{p/2}\\
    &= 1-\left(\norm {\widetilde{x}_1}^2+\norm{\widetilde{u}_1}^2-2\sum_{t=1}^{a}\widetilde{x}_{1t}\widetilde{u}_{1t}\right)^{p/2}-\left(\norm {\widetilde{x}_2}^2+\norm{\widetilde{u}_2}^2-2\sum_{t=1}^{b}\widetilde{x}_{2t}\widetilde{u}_{2t}\right)^{p/2}\\
    &= 1-\widetilde{f_u}_1(x)-\widetilde{f_u}_2(x)
    \end{align*}
Utilizing multi-index notation, we can expand
\[\widetilde{f_u}_1(x)=\sum_{\substack{\varepsilon,g\\\varepsilon+\gamma\le p/2}}(-2)^{\gamma}{\gamma\choose g}{p/2\choose \varepsilon,\gamma,p/2-\varepsilon-\gamma}\widetilde{u}_1^g\norm{\widetilde{u}_1}^{p-2\varepsilon-2\gamma}\widetilde{x}_1^g\norm{\widetilde{x}_1}^{2\varepsilon}.\]
The sum is taken over all non-negative integers $\varepsilon$ and non-negative integer vectors $g$ with $a$ entries such that $\varepsilon+\gamma\le p/2$, where $\gamma$ is the sum of the components of $g$. Similarly, we can expand
\[\widetilde{f_u}_2(x)=\sum_{\substack{\varepsilon,g\\\varepsilon+\gamma\le p/2}}(-2)^{\gamma}{\gamma \choose g}{p/2\choose \varepsilon,\gamma,p/2-\varepsilon-\gamma}\widetilde{u}_2^g\norm{\widetilde{u}_2}^{p-2\varepsilon-2\gamma}\widetilde{x}_2^g\norm{\widetilde{x}_2}^{2\varepsilon}.\]
We want to count the number of monomials in each polynomial. Since $\widetilde{f_u}_1$ is a polynomial in $\widetilde{x}_{11},\ldots,\widetilde{x}_{1a}$ and $\widetilde{f_u}_2$ is a polynomial in $\widetilde{x}_{21},\ldots,\widetilde{x}_{2b}$, the two sets of monomials are disjoint, except for the constant monomial. Let us consider $\widetilde{f_u}_1$ first. By choosing 
$\varepsilon=0$, we must count all the monomials with degree at most $p/2$. There are ${a+p/2\choose a}$ of these. If the degree is greater than $p/2$, say $p/2+c$, we only need to count the monomials formed when $\varepsilon=c$ and $\gamma=p/2-c$. Hence, the total number of monomials in $\widetilde{f_u}_1$ is
\[{a+p/2\choose a}+{a+p/2-2\choose a-1}+{a+p/2-3\choose a-1}+\cdots+{a-1\choose a-1}={a+p/2\choose a}+{a+p/2-1\choose a}.\]
Similarly, there are ${b+p/2\choose b}+{b+p/2-1\choose b}$ monomials in $\widetilde{f_u}_2$. In total, $f_u$ has ${a+p/2\choose a}+{a+p/2-1\choose a}+{b+p/2\choose b}+{b+p/2-1\choose b}-1$ monomials (we subtract $1$ as the constant monomial is counted twice). This gives an upper bound on $m$.

By finding a larger linearly independent set of polynomials, a trick first used by Blokhuis \cite{blokhuistwo}, we can lower this bound to $\binom{a+p/2}{a}+\binom{b+p/2}{b}$. We prove that the set of polynomials
\[\{f_u,\widetilde{x}_1^{m},\widetilde{x}_2^{n},1: u\in S, 0<\mu<p/2,0<\nu<p/2\}\]
where $m,n$ are integer vectors with $a,b$ entries and sums of components $\mu, \nu$ respectively, is linearly independent. The details are very similar to those in Blokhuis's \cite{Blokhuis} bound for the $s$-distance set in $\R^n$.

Suppose we have a dependence
\begin{equation}\label{eq:bigdependence}
\sum_{u\in S}a_uf_u(x)+\sum_{0<\mu<p/2}a_m\widetilde{x}_1^{m}+\sum_{0<\nu<p/2}a_n\widetilde{x}_2^{n}+\delta=0.
\end{equation}
The main claim is the following lemma.
\begin{lem}\label{lem:bigextract}
For all $m$ with $\mu<p/2$,
\[\sum_{u\in S}a_u\widetilde{u}_1^{m}=0.\]
Similarly, for all $n$ with $\nu<p/2$,
\[\sum_{u\in S}a_u\widetilde{u}_2^n=0.\]
\end{lem}
\begin{proof}
We only prove the first statement. The argument for the second statement is identical.

Suppose $\mu=t<p/2$. Consider the part of the left hand side of (\ref{eq:bigdependence}) that is homogeneous in $\widetilde{x}_{11},\ldots,\widetilde{x}_{1a}$ with degree $p-t$. Note that the monomials $\widetilde{x}_1^{m}$ do not contribute to this, so we only have to look at the part from the $f_u$. Using our expansion of $f_u$ above, the part of $f_u$ homogeneous in the $a$ variables with degree $p-t$ is
\[-\sum_{\substack{\varepsilon,g\\ 2\varepsilon+\gamma=p-t\\ \varepsilon+\gamma\le p/2}}(-2)^{\gamma}{\gamma\choose g}{p/2\choose \varepsilon,\gamma,p/2-\varepsilon-\gamma}\widetilde{u}_1^{g}\norm{\widetilde{u}_1}^{p-2\varepsilon-2\gamma}\widetilde{x}_1^{g}\norm{\widetilde{x}_1}^{2\varepsilon}.\]
The left hand side of (\ref{eq:bigdependence}) is a polynomial in the $x_i$ which is identically zero. Thus, we have
\[\sum_{\substack{\varepsilon,g\\ 2\varepsilon+\gamma=p-t\\\varepsilon+\gamma\le p/2}}(-2)^{\gamma}{\gamma\choose g}{p/2\choose \varepsilon,\gamma,p/2-\varepsilon-\gamma}\sum_{u\in S}a_u\widetilde{u}_1^{g}\norm{\widetilde{u}_1}^{p-2\varepsilon-2\gamma}\widetilde{x}_1^{g}\norm{\widetilde{x}_1}^{2\varepsilon}=0.\]
Now, substitute $x=v$, multiply by $a_v\norm{\widetilde{v}_1}^{2t-p}$, and sum over all $v\in S$.
\[\sum_{\substack{\varepsilon,g\\ 2\varepsilon+\gamma=p-t\\ \varepsilon+\gamma\le p/2}}(-2)^{\gamma}{\gamma\choose g}{p/2\choose \varepsilon,\gamma,p/2-\varepsilon-\gamma}\left(\sum_{u\in S}a_u\norm{\widetilde{u}_1}^{p-2\varepsilon-2\gamma}\widetilde{u}_1^{g}\right)^2=0.\]
This is a sum of squares where all coefficients have the same sign. So,
\[\sum_{u\in S}a_u\norm{\widetilde{u}_1}^{p-2\varepsilon-2\gamma}\widetilde{u}_1^g=0.\]
Plugging in $\gamma=t$ proves the lemma.
\end{proof}

Now, plug in $u$ into (\ref{eq:bigdependence}), multiply by $a_u$ and sum over all $u\in S$,
\[\sum_{u\in S}a_u^2+\sum_{0<\mu<p/2}a_m\sum_{u\in S}a_u\widetilde{u}_1^{m}+\sum_{0<\nu<p/2}a_n\sum_{u\in S}a_u\widetilde{u}_2^{n}+\delta\sum_{u\in S}a_u=0.\]
Applying Lemma~\ref{lem:bigextract}, we obtain $a_u=0$ for all $u$. We already know that the other polynomials in our set are linearly independent, so all coefficients vanish, as desired. This implies $m\le {a+p/2\choose a}+{b+p/2\choose b}$.
\end{proof}

\begin{rem}
This proof can be easily extended to the $\ell_p$ sum of $n$ Euclidean spaces. If we consider the product space $\E^{a_1}\oplus_p\cdots\oplus_p\E^{a_n}$, our bound is just ${a_1+p/2\choose a_1}+\cdots+{a_n+p/2\choose a_n}$.
\end{rem}

\section{Bound on equilateral sets in $\ell_p$ product spaces for $1\le p<\infty$}\label{sec:proof5} 

\begin{proof}[Proof of Theorem~\ref{thm:result5}]
Let $\{p_1,\ldots, p_m\}$ be an equilateral set in $\mathbb{E}^{a_1}\oplus_p \cdots \oplus_p \mathbb{E}^{a_n}$, and write $p_i=(\widetilde{p_i}_1,\ldots, \widetilde{p_i}_n)$ where $\widetilde{p_i}_k\in \mathbb{E}^{a_k}$ for $k=1,\ldots,n$.
 
Let $B(p)$ be the constant from Lemma~\ref{lem:jackson}. Take $c=\max(B(p),(2^{1/p}-1)^{-p})$ and set $d$ to be a positive integer satisfying \[cn\sqrt{m}<d^p<2cn\sqrt{m}.\] By Lemma~\ref{lem:jackson}, there exists an even polynomial $P$ with $P(0)=0$ and degree at most $d$ that approximates $|x|^p$. Thus, for $i\neq j$ and $k= 1,\ldots,n$, \[\left|P(\norm{\widetilde{p_i}_k-\widetilde{p_j}_k})-\norm{\widetilde{p_i}_k-\widetilde{p_j}_k}^p\right|\le \frac{B(p)}{d^p}.\]
Next, for $i=1,\ldots,m$ define the functions $f_i:\R^{a_1+\cdots+a_n}\to\R$ by
\[f_i(x)=1-\sum_{k=1}^n P(\norm{\widetilde{p_i}_k-\widetilde{x}_k}).\]
Let $M$ be the $m\times m$ matrix given by $m_{ij}=f_i(p_j)$. Since $P(0)=0$, we have $m_{ii}=1$. For $i\neq j$, \begin{align*}
    |m_{ij}| &= \left|\sum_{k=1}^n \norm{\widetilde{p_i}_k-\widetilde{p_j}_k}^p-\sum_{k=1}^n P(\norm{\widetilde{p_i}_k-\widetilde{p_j}_k})\right|\\
    &\le \sum_{k=1}^n \left | \norm{\widetilde{p_i}_k-\widetilde{p_j}_k}^p-P(\norm{\widetilde{p_i}_k-\widetilde{p_j}_k})\right |\\
   &\le \frac{nB(p)}{d^p}\\
  &< \frac{1}{\sqrt{m}},
\end{align*} with the last inequality following from our conditions on $c$ and $d$. Because $M$ is symmetric, we can apply Lemma~\ref{lem:ranklemma} to get $\rank M\ge m/2$.

Now, we look for an upper bound on the rank. We can write any point $x\in \E^{a_1}\oplus_p\cdots\oplus_p\E^{a_n}$ as
\[x=(\widetilde{x}_1^{(1)},\widetilde{x}_1^{(2)},\ldots, \widetilde{x}_1^{(a_1)},\widetilde{x}_2^{(1)},\widetilde{x}_2^{(2)},\ldots, \widetilde{x}_2^{(a_2)},\ldots, \widetilde{x}_n^{(1)},\widetilde{x}_n^{(2)},\ldots, \widetilde{x}_n^{(a_n)}).\]
Because $P$ is even, all terms in the expansion of $P(\norm{\widetilde{p_i}_k-\widetilde{x}_k})$ have integer exponent, i.e., $f_i$ is actually a polynomial. Additionally, since $P$ has degree at most $d$, each $f_i$ is in the span of the set \[\mathcal{S}=\{1, \sum_{k=1}^n \norm{\widetilde{x}_k}^d\}\cup \mathcal{S}_{1}\cup \mathcal{S}_{2}\cup \cdots\cup \mathcal{S}_{n},\]
where each set $\mathcal{S}_{i}$ consists of all the monomials with degree less than $d$ formed by $\widetilde{x}_i^{(1)}, \ldots, \widetilde{x}_i^{(a_i)}$.
Thus, the $f_i$ are spanned by at most $2+\sum_{k=1}^n \binom{a_k+d-1}{a_k}-n$ polynomials. Because every row vector of $M$, each of the form $(f_i(p_1),\ldots, f_i(p_m))$, belongs to the subspace spanned by the set $\{(f(p_1),\ldots,f(p_m)):f\in \mathcal{S}\}$, we have
\begin{align*}
\rank M &\le 2+\sum_{k=1}^n\binom{a_k+d-1}{a_k}-n\\
&\le \sum_{k=1}^n \binom{a_k+d-1}{a_k}\\
&\le n\binom{a+d-1}{a}\\
&< ne^a\left(1+\frac{d}{a}\right)^a,
\end{align*}
where we let $\displaystyle a=\max_{1\le i\le n}a_i$ and use the well-known bound ${n\choose k}\le (en/k)^k$ for $1\le k\le n$. Now, recall that $d\ge B(p)^{1/p}\ge p$ and $2p>a$, so $d>a/2$. Thus, \[\frac{m}{2}\le\rank M< n\left(\frac{3ed}{a}\right)^a.\] Combining this inequality with the condition $d^p<2cn\sqrt{m}$ yields \[m<c_{p,a}n^{\frac{2p+2a}{2p-a}},\]
as desired.
\end{proof}

\section{Acknowledgements}
The authors would like to thank the MIT PRIMES program for making this research possible. The authors would also like to thank Larry Guth for introducing them to the problem, as well as Yufei Zhao for some helpful suggestions. 

\bibliographystyle{abbrv}
\bibliography{bibfile}

\end{document}